\documentclass[a4paper,reqno,oneside]{article}

\usepackage{amsthm,amsmath,amssymb,enumitem,mathrsfs}
\usepackage{cancel}

\usepackage[latin1]{inputenc}

\usepackage{marginnote}
\usepackage{inconsolata}

\usepackage[notref,notcite,final]{showkeys}

\usepackage[dvipsnames]{xcolor}
\usepackage[draft, markup=sout, authormarkup=brackets]{changes}

\usepackage[colorlinks=true,citecolor=magenta,linkcolor=magenta]{hyperref}

\usepackage[modulo]{lineno}
\modulolinenumbers[1]

\newcommand*\patchAmsMathEnvironmentForLineno[1]{%
	\expandafter\let\csname old#1\expandafter\endcsname\csname #1\endcsname
	\expandafter\let\csname oldend#1\expandafter\endcsname\csname end#1\endcsname
	\renewenvironment{#1}%
	{\linenomath\csname old#1\endcsname}%
	{\csname oldend#1\endcsname\endlinenomath}}%
\newcommand*\patchBothAmsMathEnvironmentsForLineno[1]{%
	\patchAmsMathEnvironmentForLineno{#1}%
	\patchAmsMathEnvironmentForLineno{#1*}}%
\AtBeginDocument{%
	\patchBothAmsMathEnvironmentsForLineno{equation}%
	\patchBothAmsMathEnvironmentsForLineno{align}%
	\patchBothAmsMathEnvironmentsForLineno{flalign}%
	\patchBothAmsMathEnvironmentsForLineno{alignat}%
	\patchBothAmsMathEnvironmentsForLineno{gather}%
	\patchBothAmsMathEnvironmentsForLineno{multline}%
}

\newcommand{\eps}{\varepsilon}
\newcommand{\epsi}{\varepsilon}
\renewcommand{\epsilon}{\varepsilon}

\newcommand{\rhs}{right hand side }

\newcommand{\erfc}{{\rm erfc}}

\DeclareMathOperator{\deriv}{d}

\DeclareMathOperator{\dd}{d}

\let\TeXchi\chi
\def\chi{{\setbox0 \hbox{\mathsurround0pt
			$\TeXchi$}\hbox{\raise\dp0 \copy0 }}}

\newcommand{\scrH}{{\mathscr H}}

\newcommand{\calH}{{\mathcal H}}

\newcommand{\dit}{\deriv\!t}

\newcommand{\derivt}{\frac{\deriv\!{}}{\dit}}

\makeatletter
\def\thmheadbrackets#1#2#3{%
	\thmname{#1}\thmnumber{\@ifnotempty{#1}{ }\@upn{#2}}%
	\thmnote{ {\the\thm@notefont[#3]}}}
\makeatother

\newtheoremstyle{brakets}
{}
{}
{\itshape}
{}
{\bfseries}
{.}
{ }
{\thmheadbrackets{#1}{#2}{#3}}

\theoremstyle{brakets}

\newtheorem{thm}{Theorem}[]
\newtheorem*{thm-nn}{Theorem}
\newtheorem{lem}[thm]{Lemma}
\newtheorem*{lem-nn}{Lemma}

\newtheorem{cor}[thm]{Corollary}
\newtheorem*{cor-nn}{Corollary}
\newtheorem*{remark}{Remark}

\newtheorem*{ex-nn}{Example}

\newcommand{\grad}{\nabla_x} 
\newcommand{\lap}{\Delta_x} 

\newcommand{\R}{\mathbb{R}}

\newcommand{\sprod}[2]{\langle #1 , #2 \rangle} 

\DeclareMathOperator{\ds}{ds}

\newcommand{\ints}[1]{\int_{\T^d}#1 \,dx} 

\newcommand{\intsv}[1]{\int_{\T^d}\int_{\R^d} #1\, dv\,dx }

\newcommand{\intv}[1]{\int_{\R^d} #1\, dv}

\renewcommand{\phi}{\varphi}

\newcommand{\bfQ}{\mathbf{Q}}

\def\sfT{\mathsf{T}}
\def\sfL{\mathsf{L}}

\def\sfA{\mathsf{A}}

\newcommand{\finf}{f_\infty}
\newcommand{\rinf}{\rho_\infty}
\newcommand{\Tinf}{T_\infty}

\newcommand{\T}{\mathbb{T}}
\newcommand{\hf}{\hat{f}}
\newcommand{\hT}{\hat{T}}
\newcommand{\hr}{\hat{\rho}}
\newcommand{\tf}{\tilde{f}}
\newcommand{\tT}{\tilde{T}}
\newcommand{\tr}{\tilde{\rho}}

\newcommand{\nn}{\nonumber}

\allowdisplaybreaks

\usepackage{empheq}
\usepackage[most]{tcolorbox}
\newtcbox{\mymath}[1][]{%
	nobeforeafter, math upper, tcbox raise base,
	enhanced, colframe=black!50,
	colback=black!10, boxrule=1pt,
	#1}

\usepackage[affil-it]{authblk}

\begin{document}
	
	
\title{\vspace{-2cm}\textbf{Thermalization of a rarefied gas with total energy conservation: 
		existence, hypocoercivity, macroscopic limit}}

\author[1]{Gianluca Favre\footnote{E-mail: {\tt gianluca.favre@univie.ac.at}}}

\author[1]{Christian Schmeiser\footnote{E-mail: {\tt christian.schmeiser@univie.ac.at}}}

\author[2]{Marlies Pirner\footnote{E-mail: {\tt marlies.pirner@mathematik.uni-wuerzburg.de}}}

\affil[1]{Faculty of Mathematics, University of Vienna,	Oskar-Morgenstern-Platz 1, 1090 Wien, Austria}

\affil[2]{Department of mathematics, W{\"u}rzburg University, Emil Fischer Str. 40, 97074 W{\"u}rzburg, Germany}


		
\date{}
\maketitle	
	
\begin{abstract}
\noindent The thermalization of a gas towards a Maxwellian velocity distribution with the background temperature is described by a kinetic relaxation model. The sum 
of the kinetic energy of the gas and the thermal energy of the background are conserved, and the heat flow in the background is governed by the Fourier law.

For the coupled nonlinear system of the kinetic and the heat equation, existence of solutions is proved on the one-dimensional torus. Spectral stability 
of the equilibrium is shown on the torus in arbitrary dimensions by hypocoercivity methods. The macroscopic limit towards a nonlinear cross-diffusion
problem is carried out formally.
\end{abstract}	\medskip


\vspace{2mm}

\noindent {\bf AMS (MOS) subject clas\-si\-fi\-ca\-tion:}%
~~ 82C40, 35M30, 35Q70, 35A01, 35B40. 
%
	
\section{Introduction}
We consider a gas in a periodic box exchanging energy with a background, where heat conduction is governed by the Fourier law. The energy exchange is the 
consequence of thermalizing scattering events, where post-collisional velocities are sampled from a Maxwellian distribution with zero mean (assuming that the 
background is at rest) and with the background temperature. Collisions between gas molecules are neglected. A motivation for this work is a first small step towards
the extension of kinetic transport models for chemical reaction networks (see, e.g., \cite{BisDes,FavreSchm}) to exothermic or endothermic reactions where the 
energy balance needs to be considered.

The mathematical model consists of a BGK type kinetic equation for the gas with relaxation towards the above described Maxwellian, and an inhomogeneous heat equation for the background temperature, where the inhomogeneity is chosen to ensure conservation of the total energy, i.e. the sum of the kinetic energy of the
gas and the thermal energy of the background. 


Our first main result is global existence of a solution of the initial value problem in the one-dimensional case. The latter restriction results from the fact that 
boundedness of the background temperature is needed for controlling the nonlinearities. However, by the energy transfer the kinetic energy density
of the gas, which is only integrable, appears as an inhomogeneity in the heat equation, producing a bounded temperature only in dimension one.
We are not aware of similar results in the literature. However, for macroscopic models of similar physical settings see \cite{Mar-Hitt-Hask-Miel}, where an
existence result also requires considerable effort.

The second main result is spectral stability of the global equilibrium in arbitrary dimensions. The first step towards this result is to understand the thermodynamic
structure of the problem. This is facilitated by an interpretation of the problem as an approximation for a two-species gas mixture in the limit of disparate 
molecule cross sections. In the limit, the large particle species becomes the equilibrium background, and the standard entropy/entropy dissipation structure carries
over to the limit problem. Entropy dissipation works towards a constant background temperature and a gas with (Maxwellian) equilibrium velocity distribution and
arbitrary macroscopic density. For a global equilibrium the macroscopic density is also expected to be constant as a result of mixing by the kinetic transport.
The global equilibrium can be determined uniquely from mass conservation in the gas and from total energy conservation. Because of the non-definiteness of the
entropy dissipation, spectral stability is a {\em hypocoercivity} result. It is shown by employing the abstract $L^2$-hypocoercivity approach of 
\cite{DolMouSchmHypoMass}, which relies on the construction of a Lyapunov functional by augmenting the entropy of the linearized problem such that also
the kinetic transport contributes to the decay.

The macroscopic limit for dominating gas-background scattering is carried out formally. It is of diffusive nature and takes the form of a nonlinear cross-diffusion 
system for the macroscopic gas density and the background temperature. An existence result for this limiting problem is the subject of the parallel effort \cite{FJSZ-2020}.


It is an essential difference to the Boltzmann BGK model (see \cite{PERTHAME1989191,perthame1993weighted} for existence proofs) that the temperature 
in the Maxwellian is not the temperature of the gas but of the background. This explains the dimension restriction in our result. Modified BGK models have been treated in \cite{bosi2009bgk, choi2018global, yun2015classical, park2019cauchy, klingenberg2018existence}. We use some of the ideas from these contributions.

Largely motivated by the book of Villani \cite{Villani}, the literature on hypocoercivity has been growing considerably, with most of the approaches based on the
construction of suitable Lyapunov functionals. The $H^1$-based approach initiated in \cite{MouNeu} and expanded in \cite{Villani} is strongly motivated by
the theory of hypoellipticity. Recently it has been extended for certain model problems to prove sharp decay rates \cite{Achleitner2015,Arnold2014}.
The $L^2$ approach of \cite{DolMouSchmHypoMass} has been strongly motivated by Herau \cite{Herau} and is close to the Kawashima
modified function approach \cite{Kawashima}. Its abstract formulation permits applications not only to kinetic transport equations, but also to certain coupled systems
like, e.g., in \cite{DBGV16,FavreSchm,LP2019} or in this work.

%


The rest of this article is structured as follows. In Section \ref{secModel}, the model is formulated including an outline of a derivation from a model for a mixture of
two gases with disparate collision cross sections, motivating the thermodynamic structure, which is also reduced to the linearization around the global equilibrium.
In the ensuing Section \ref{sec hypo main} the main results are presented. The proofs of hypocoercivity of the linearized system and of global existence for the nonlinear system in one dimension are contained in Sections \ref{sec4} and, respectively, \ref{sec existence}.

\section{The Model}
\label{secModel}
We define the Maxwellian 
$$M(T)(v) := \big( \pi \, T\big)^{-d/2} \exp\bigg({-\frac{|v|^2}{ T}}\bigg)$$
where $v \in \R^d$ and $x\in\T^d$, the flat torus represented by $[0,1]^d$. In the following, its moments
\begin{eqnarray}
  && \intv{M(T)} = 1 \,,\qquad \intv{|v|^2 M(T)} = \frac{d}{2}T \,,\nn\\
  && \intv{|v|^4 M(T)} = \frac{d(d+2)}{4}T^2\,,\label{M-moments}
\end{eqnarray}
will be needed.
$T(x,t)$ denotes the temperature of the background and $f(x,v,t)$ the phase space number density of the gas. Its macroscopic density $\rho(x,t)$ 
is defined by $\rho(x,t) = \intv{f(x,v,t)}$. Moreover, the macroscopic energy density $E(x,t)$ is defined by $E(x,t) = \intv{|v|^2 f(x,v,t)}$.

All the variables have been nondimensionalized such that $f$ and $T$ satisfy the system
\begin{eqnarray}
\label{Kin}
\partial_t f + v \cdot \grad f &=& \rho\, M(T) - f \,,\\
\label{Heat}
\partial_t T - D \lap T &=& \intv{|v|^2 (f - \rho\, M(T))} = E - \frac{d}{2}\rho T\,, 
\end{eqnarray}
with the nondimensionalized heat conductivity $D$ of the background, subject to initial conditions
\begin{equation}\label{IC}
f(x,v,0)=f_0(x,v)\,, \qquad T(x,0)=T_0(x) \,.
\end{equation}
The first equation describes the kinetic transport of the gas and its relaxation towards the Maxwellian equilibrium with vanishing average velocity and with the
temperature of the background. The second equation describes heat conduction within the background 
as well as the exchange of energy with the gas. The right-hand side is chosen such that the total energy is locally conserved:
\begin{equation}
\label{Cont eq}
\partial_t \left(E + T \right) + \grad \cdot \left(\intv{v \, |v|^2 f} - D \grad T \right) = 0 \,.
\end{equation}
The normalization of the Maxwellian, $\intv M=1$, implies mass conservation for the gas:
\begin{equation} \label{Cont1}
\partial_t \rho + \grad \cdot \left(\intv{v f}  \right) = 0\,.
\end{equation}

\subsection*{Entropy and equilibrium}

The thermodynamic structure of \eqref{Kin}, \eqref{Heat} can be determined by interpreting it as an approximation for a more general model.
Consider a nondimensionalized kinetic model for a two-component gas mixture:
\begin{eqnarray}
\label{Boltz1}
\partial_t f + v \cdot \grad f &=& Q_{12}(f,g) + \epsilon Q_{11}(f) \\
\label{Boltz2}
\partial_t g + v \cdot \grad g &=& Q_{21}(g,f) + \epsilon^{-1} (Q_{22}(g) + Q_{el}(g))\,,
\end{eqnarray}
where $f(x,v,t)$ and $g(x,v,t)$ are the phase space distributions of the two components 1 and, respectively, 2. The small dimensionless parameter 
$\epsilon$ results from the scaling and could have the interpretation of the ratio of cross section areas of the two species. The explicit form of the 
operators on the right hand sides, acting only in the velocity direction, are not essential for our argument and therefore omitted. The coupling operators $Q_{12}$ 
and $Q_{21}$ describe binary collisions 
between a particle of component 1 and a particle of component 2. On the other hand, the operators $Q_{11}$ and $Q_{22}$ are models for binary
collisions within the components, while $Q_{el}$ describes elastic, directionally unbiased collisions with a nonmoving background. All the collision
processes are assumed to conserve the particle numbers of both components as well as the total kinetic energy, in particular
\begin{eqnarray}
    &&\intv{Q_{12}(f,g)} = 0 \,,\label{mass}\\
    &&\intv{|v|^2(Q_{12}(f,g) + Q_{21}(g,f))} =  0 \,,\label{energy1}\\
    &&\intv{|v|^2(Q_{22}(g) + Q_{el}(g))} = 0 \,.\label{energy2}
\end{eqnarray}
Finally we expect the entropy
$$
   H(f,g) = \intsv{\big(f \ln f + g \ln g \big)}
$$
to be nonincreasing along solutions of \eqref{Boltz1}, \eqref{Boltz2}. 

In the limit $\epsilon\to 0$, the distribution $g$ can be expected to satisfy $Q_{22}(g) + Q_{el}(g) = 0$, assumed to imply that $g$ is a Maxwellian 
distribution with vanishing mean velocity, the latter as a result of the collisions with the background. We choose a solution with constant macroscopic 
density equal to one, i.e. $g(x,v,t) = M(T(x,t))(v)$ in the limit $\epsilon\to 0$. We also assume that particles of the two components have the same mass, 
with the consequence that $Q_{12}(f, M(T))= 0$ implies that $f$ is of the form $f(x,v,t) = \rho(x,t)M(T(x,t))(v)$. The simplest model satisfying this and
also particle conservation \eqref{mass} is the relaxation model
$$
   Q_{12}(f,M(T)) = \rho M(T) - f \,,\qquad \mbox{with } \rho = \intv{f} \,,
$$
giving \eqref{Kin} as the limit of \eqref{Boltz1}.
An equation for $T$ is then obtained by multiplication of \eqref{Boltz2} by $|v|^2$, integration with respect to $v$, and passing to the limit
(using \eqref{energy1} and \eqref{energy2}):
\begin{eqnarray*}
  &&\partial_t \left( \intv{|v|^2 M(T)}\right) = \intv{|v|^2 Q_{21}(M(T),f)} \\
  &&= - \intv{|v|^2 Q_{12}(f,M(T))} = \intv{|v^2|(f-\rho M(T))} \,.
\end{eqnarray*}
This is not quite \eqref{Heat} since, on the one hand, the heat conduction term is missing and, on the other hand, there is an extra factor $d/2$ since,
by \eqref{M-moments}, $\intv{|v|^2 M(T)} = dT/2$.

For the entropy we obtain
$$
    H(f,M(T)) = \ints{\left(\intv{f\ln f} - \frac{d}{2}\ln T\right)} + \mbox{ const}
$$
Considering the above mentioned factor $d/2$, this suggests that 
\begin{equation}
\calH(f,T) =\ints{\left(\intv{f \ln f } - \ln T\right)} 
\label{H1}
\end{equation} 
is an entropy for \eqref{Kin}, \eqref{Heat}, since the integrand is a convex function of $T$, thus compatible with the heat conduction term,
which can be interpreted as the result of a Chapman-Enskog expansion. Indeed, computing the time derivative we obtain
\begin{align*}
\nn
\frac{d}{dt}\calH(f,T) & = \intsv{ (\rho M(T) - f) \left( \ln f + |v|^2  \frac{1}{T}  \right) } - \ints{ D\,    \frac{\lap T}{T}  } \\
   & = - \intsv{ (f - \rho M(T)) \ln \frac{f}{\rho M(T)} } - D \ints{\frac{|\grad T|^2}{T^2}} \le 0\,.
\end{align*}
As expected, the heat conduction produces a contribution to the entropy dissipation, whose form suggests convergence of solutions of \eqref{Kin},
\eqref{Heat} as $t\to\infty$ to an equilibrium 
$$
   (\finf,\Tinf) = (\rho_\infty M(\Tinf), \Tinf) \,,
$$ 
with constant $\Tinf$ and possibly position
dependent $\rho_\infty$. However, considering the transport term in \eqref{Kin}, which does not contribute to the entropy dissipation, for an equilibrium
also $\rho_\infty$ has to be constant. The values of the two constants can be determined from the initial data \eqref{IC} by integration of \eqref{Cont eq} 
and \eqref{Cont1} with respect to $x$, i.e. by conservation of the total mass of the gas and by conservation of the total energy:
$$
   \rho_\infty |\Omega| = \intsv{f_0} \,,\qquad \left( \rho_\infty \frac{d}{2} + 1\right) \Tinf |\Omega| = \ints{ \left(\intv{|v|^2 f_0} + T_0\right)} \,.
$$
The relative entropy (relative to the equilibrium) is given by
\begin{equation}
\label{Rel entropy}
{\calH}\bigl((f,T)|(\finf,T_\infty)\bigr) = \ints{\left( \intv{f \ln \frac{f}{\finf} } - \ln \frac{T}{\Tinf} + \frac{T}{\Tinf} - 1 \right) } \,.
\end{equation}
Since it differs from $\calH(f,T)$ only by adding a constant and a constant multiple of the total energy, it has the same dissipation.

\subsection*{Linearization}

The perturbations $\hf := f - \finf$, $\hT := T - \Tinf$ satisfy 
\begin{equation}\label{cons-zero}
  \intsv{\hf} = \ints{\hr} = 0 \,,\qquad \ints{\bigg(\intv{|v|^2 \hf} + \hT \bigg)} = 0 \,.
\end{equation}
Linearization of \eqref{Kin}, \eqref{Heat} around the equilibrium gives
\begin{align}
 \partial_t \hat{f} +  v \cdot \grad \hat{f} &= \bfQ_L(\hf,\hT) \,, \label{Kinlin}\\ 
\partial_t \hat{T} - D~\lap \hat{T}  &= - \intv{|v|^2 \bfQ_L(\hf,\hT)} \,,\label{Heatlin}
\end{align}
with the linearized collision operator 
\begin{equation*}
\bfQ_L(\hf,\hT) = \hr M(\Tinf) + \finf \frac{\hT}{\Tinf} \bigg( \frac{|v|^2}{\Tinf} - \frac{d}{2} \bigg) - \hf \,,
\end{equation*}
which shares the mass conservation property $\intv{\bfQ_L(\hf,\hT)} = 0$ with its nonlinear counterpart.
An entropy for the linearized system is obtained as the quadratic approximation of the relative entropy \eqref{Rel entropy} close to equilibrium: 
\begin{equation}
\label{quadr entr}
\calH_L(\hf,\hT) = \ints{ \bigg( \intv{ \frac{\hf^2}{2\finf} } + \frac{\hT^2}{2\Tinf^2} \bigg) } \,.
\end{equation}
Its dissipation is computed as 
\begin{eqnarray}
  \frac{d}{dt} \calH_L(\hf,\hT) &=& \intsv{ \frac{\hf}{\finf} \bfQ_L(\hf,\hT)} \nn\\
  && + \ints{\frac{\hT}{\Tinf^2} \left( D\Delta_x \hT - \intv(|v|^2\bfQ_L(\hf,\hT))\right)} \nn\\
  &=& - \intsv{\frac{\bfQ_L(\hf,\hT)^2}{\finf}} - \frac{D}{\Tinf^2} \ints{|\nabla_x \hT|^2} \label{lin-entr-diss}\\
  && + \intsv{\bfQ_L(\hf,\hT) \left( \frac{\hr}{\rinf} - \frac{d\,\hT}{2\Tinf}\right)} \,,\nn
\end{eqnarray}
where the last line vanishes by the mass conservation property of $\bfQ_L$. As a plausibility check, the dissipation \eqref{lin-entr-diss} 
can also be obtained as the quadratic approximation close to equilibrium of the nonlinear entropy dissipation.

\section{Main results}
\subsection*{Hypocoercivity for the linearized problem} \label{sec hypo main}

As for the nonlinear problem, the entropy dissipation \eqref{lin-entr-diss} is not definite. It vanishes whenever $\bfQ_L(\hf,\hT) = 0$, i.e.
$\hf = \hr M(\Tinf) + \finf \frac{\hT}{\Tinf} \left( \frac{|v|^2}{\Tinf} - \frac{d}{2} \right)$, and $\hT$ is independent of $x$. Using this form in the conservation
equations \eqref{cons-zero} implies $\hT=0$. However, $\hr$ can be an arbitrary function of $x$ with average zero. Despite this nondefiniteness we expect 
hypocoercivity, i.e. exponential convergence of $(\hf,\hT)$ to zero as $t\to\infty$.
The following hypocoercivity result will be proved in Section \ref{sec4} in the natural functional analytic setting of the weighted $L^2$ space induced 
by the quadratic entropy $\calH_L$.

\begin{thm} {\bf (Hypocoercivity)}
	\label{Thm Hypo}
Let $D>0$. Then there exist positive constants $C$ 
and $ \lambda$ such that solutions of \eqref{cons-zero}, \eqref{Kinlin},\eqref{Heatlin} satisfy
\begin{equation*}
 \calH_L(\hf,\hT) \le C e^{-\lambda t} \calH_L(\hf(t=0),\hT(t=0)) \,.
\end{equation*}
\end{thm}

\subsection*{Formal macroscopic limit: cross diffusion}

In this section, we derive the formal macroscopic limit under the diffusive macroscopic scaling $x \to x\, \epsi^{-1}$, $t \to t \, \epsi^{-2} $ with $0 < \epsi \ll 1$ . 
 The rescaled version of \eqref{Kin}, \eqref{Heat} reads
\begin{align}
\label{Kineps}
\epsilon^2 \partial_t f_{\varepsilon} + \epsilon v \cdot \grad f_{\varepsilon} &= \rho_{\varepsilon} M(T_{\varepsilon}) - f_{\varepsilon} \\
\label{Heateps}
\epsi^2 \big(\partial_t T_{\varepsilon} - D~\lap T_{\varepsilon} \big) &= \intv{|v|^2 \big(f_{\varepsilon} - \rho_{\varepsilon} M(T_{\varepsilon}) \big)} \,.
\end{align}
The formal limit $\varepsilon \rightarrow 0$ leads to 
\begin{align*}
f_{\varepsilon}(x,v,t) \rightarrow f_0(x,v,t)= \rho_0(x,t) M(T_0(x,t))(v) \,.
\end{align*}
Now, we split $f_{\varepsilon} = \rho_{\varepsilon} M(T_{\varepsilon}) + \varepsilon R_{\varepsilon}$ and insert this into \eqref{Kineps}, \eqref{Heateps}, we obtain 
\begin{align}
\varepsilon \partial_t f_{\varepsilon} +  v \cdot \nabla_x f_{\varepsilon} &= - R_{\varepsilon} \,,\label{Kin/eps}\\
\varepsilon \partial_t T_{\varepsilon} - \varepsilon D \lap T_{\varepsilon} &= - \intv{  |v|^2 R_{\varepsilon} } \,.\label{Heat/eps} 
\end{align}
The formal limit of \eqref{Kin/eps}, 
\begin{align}\label{R0}
R_0 = -v \cdot \nabla_x (\rho_0 M(T_0)) \,,
\end{align}
is compatible with the limit of \eqref{Heat/eps}, since $R_0$ is an odd function of $v$.

Equations for $\rho_0$ and $T_0$ will be derived as the limits of the conservation laws \eqref{Cont eq} and \eqref{Cont1}.
We integrate  \eqref{Kineps} with respect to $v$ and divide it by $\varepsilon^2$:
$$
\partial_t \rho_\eps + \nabla_x\cdot\intv{ v R_\eps } =0 \,,
$$
where we have used $\intv{ R_{\varepsilon} }=0$ and $\intv{ v \rho_{\varepsilon} M_{\varepsilon} } =0$. With \eqref{R0} we obtain in the limit $\varepsilon \rightarrow 0$ the macroscopic equation
\begin{equation}\label{mass0}
\partial_t \rho_0 - \frac{1}{2}\lap (\rho_0 T_0) =0 \,.
\end{equation}
Similarly, the energy conservation equation can be written as 
$$
   \partial_t \left( \intv{|v|^2f_\eps} + T_\eps \right) + \nabla_x\cdot \left( \intv{v|v|^2 R_\eps} - D\nabla_x T_\eps\right) = 0 \,.
$$
The limit $\eps\to 0$ gives
\begin{equation}\label{energy0}
  \partial_t \left( \left( \frac{\rho_0 d}{2} + 1\right)T_0\right) - \lap \left(\left( \frac{3d}{4}\rho_0 T_0 + D\right)T_0\right) = 0
\end{equation}
The nonlinear cross diffusion system \eqref{mass0}, \eqref{energy0} has been studied in \cite{FJSZ-2020}.


\subsection*{Global existence for the nonlinear problem in one dimension}

The main issue in proving an existence result is the control of the nonlinear coupling terms. 

A mild formulation will be used for the initial value problem for the temperature equation \eqref{Heat}:
\begin{equation}\label{mild-heat}
    T(x,t) = \int_0^1 G(x-y,t)T_0(y)dy + \int_0^t \int_0^1 G(x-y,t-s)\left(E - \frac{d}{2}\rho T\right)(y,s)dy\,ds \,,
\end{equation}
with Green's function for the heat operator on the torus $\T^1$, given by
\begin{equation}\label{Green}
   G(x,t) = (4D\pi t)^{-1/2} \sum_{k\in\mathbb{Z}} \exp\left(-\frac{(x+k)^2}{4Dt} \right) \,.
\end{equation}
The setting of Corollary \ref{cor:heat} will be used with $E,\rho\in L^\infty_t L^1_x$. This does not provide a positive lower bound for $T$. Therefore
the following distributional formulation of the term $\rho M(T)$ in the kinetic equation \eqref{Kin} will be used:
\begin{equation}\label{rhoM-distr}
  \langle \rho M(T),\phi\rangle = \int_{\mathbb{T}^1\times\R\times (0,\infty)} \rho(x,t)M(1)(w)\phi(x,w\sqrt{T},t)\dd x \dd w \dd t
\end{equation}

\begin{thm} \label{Thm_Ex}
Let $\tau>0$ and $\mathbb{T}^1$ the one-dimensional torus. Let the initial data satisfy $T_0\in C^1_+(\mathbb{T}^1)$
and $(1+|v|^q + (\log f_0)_+)f_0 \in L^1_+(\mathbb{T}^1\times\R)$ for some $q>2$.
Then there exists a solution $(f,T)$ of \eqref{Kin}--\eqref{IC}, where 
\begin{itemize}
\item $(1+|v|^2+\log f)f\in L^\infty((0,\tau); L^1(\mathbb{T}^1\times\R))$, $T\in C^{1/2}(\mathbb{T}^1\times [0,\tau])$,
\item the Cauchy problem for the kinetic equation \eqref{Kin} is interpreted in the distributional sense with the interpretation \eqref{rhoM-distr} for the term $\rho M(T)$,
\item the Cauchy problem for the heat equation \eqref{Heat} is interpreted in the mild sense \eqref{mild-heat}.
\end{itemize}
\end{thm}

\noindent The proof is given in Section \ref{sec existence}. 


%

\section{Proof of Hypocoercivity}
\label{sec4}

In this section, we employ the abstract hypocoercivity approach of \cite{DolMouSchmHypoMass} to prove Theorem \ref{Thm Hypo}, i.e. exponential decay 
for solutions of the linearised system \eqref{cons-zero}--\eqref{Heatlin}.

We write \eqref{Kinlin}, \eqref{Heatlin} in the abstract form
\begin{equation}
\label{abstract formulation}
\partial_t F + \sfT F = \sfL F \,,
\end{equation}
with $F = (\hf,\hT)$, where we separate the dissipative \textit{collision operator} $\sfL$ from the  \textit{transport operator}
$\sfT$ with the definitions
\begin{equation*}
\sfT F =
\begin{pmatrix}
v \cdot \grad \hf \\ 0
\end{pmatrix} \,, \qquad
\sfL F=
\begin{pmatrix}
\bfQ_L(F) \\ 
D \, \lap \hT - \intv{|v|^2 \bfQ_L(F)}
\end{pmatrix} \,.
\end{equation*}
The form of the quadratic entropy \eqref{quadr entr} suggests the introduction of the weighted scalar product
\begin{equation}
   \sprod{F}{G} = \ints{ \bigg(\intv{ \frac{\hf \hat{g}}{\finf} } + \frac{\hT\hat\theta}{\Tinf^2}\bigg) } \,,
\label{scalar}
\end{equation}
where $F=(\hf,\hT)$ and $G=(\hat{g},\hat\theta)$, and of the induced norm $\|\cdot\|$. We shall also use
$$
    \|\hf\|_v^2 := \intv{\frac{\hf^2}{\finf}} \,,\qquad \|\cdot\|_x := \|\cdot\|_{L^2(\Omega)} \,,\qquad \|\hf\|_{v,x} := \left\|\|\hf\|_v\right\|_x
$$ 
such that 
$$
    \calH_L(F) = \|F\|^2 =\|\hf\|_{v,x}^2 + \frac{\|\hT\|_x^2}{\Tinf^2} \,.
$$
The Hilbert space of all $F=(\hf,\hT)$ satisfying $\|F\|<\infty$
and \eqref{cons-zero} will be denoted by $\mathbb{H}$. Analogously to the derivation of the entropy dissipation \eqref{lin-entr-diss} it is easily
shown that $\sfL$ is symmetric and $\sfT$ is antisymetric with respect to $\sprod{\cdot}{\cdot}$.

We define the orthogonal projector $\Pi$ on the set of \emph{local equilibria} satisfying $\sfL F = 0$.
With $\Pi F = (\tf,\tT)$ we have $\bfQ_L(\tf, \tT) =0$ and $\lap \tT = 0$, the former implying
\begin{equation*}
\tf = \tr M(\Tinf) + \finf \frac{\tT}{\Tinf} \bigg( \frac{|v|^2}{\Tinf} - \frac{d}{2} \bigg) \,,
\end{equation*}
and the latter that $\tT$ is constant. Since $(\tf,\tT)$ has to satisfy \eqref{cons-zero}, it follows that $\ints{\tr} = 0$ and $\tT=0$. Finally the orthogonality 
of the projection implies $\tr = \hr$ and thus,
$$
    \Pi F = (\hr M(T_\infty), 0) \,.
$$
The entropy dissipation 
\begin{equation}
\label{def LFF}
\frac{d}{dt} {\calH}_L = \sprod{\sfL F}{F}\,,
\end{equation}
seems to indicate convergence to a local equilibrium, but not to the global equilibrium zero. In order to show convergence to the global equilibrium, we  define the modified entropy $\scrH[F]$ as done in \cite{DolMouSchmHypoMass},
\begin{equation*}
\scrH[F]=||F||^2 + \delta \sprod{\sfA F}{F} \,,
\end{equation*}
with
\begin{equation*}
\sfA = \big[ \mathbf{I} + (\sfT \Pi)^*(\sfT \Pi) \big]^{-1} (\sfT \Pi)^* \,,
\end{equation*}
and $\delta > 0$ to be determined later.
In \cite{DolMouSchmHypoMass} it has been shown that the operator norm of $\sfA$ is bounded by $1/2$, such that
\begin{equation}
\label{equival}
\frac{1-\delta}{2} \|F\|^2 \le \scrH[F] \le \frac{1+\delta}{2} \|F\|^2 \,,
\end{equation}
i.e. $\scrH[\cdot]$ is equivalent to $\|\cdot\|^2$ for $\delta < 1$.
The time derivative of the modified entropy along solutions of \eqref{abstract formulation} can be computed as
\begin{equation}
\label{derivt entr}
\frac{d}{dt} \scrH[F] = \sprod{\sfL F}{F} - \delta \sprod{\sfA \sfT \Pi F}{F} + \delta \sprod{\sfT \sfA F}{F} - \delta \sprod{\sfA \sfT (\mathbf{I}-\Pi)F}{F} + \delta \sprod{\sfA \sfL F}{F} \,.
\end{equation}
The idea is that the first two terms on the \rhs provide the desired coercivity, whereas the remaining three terms are perturbations to be controlled by the first two. In the following auxiliary results the assumptions used in the abstract approach of \cite{DolMouSchmHypoMass} are verified. We shall make use of
the \emph{Poincar\'e inequality} on $\Omega$: There exists a constant $c_P>0$, such that
\begin{equation}\label{Poincare}
   \|\nabla_x \theta\|_x^2 \ge c_P \left\| \theta - \frac{1}{|\Omega|}\ints{\theta}\right\|_x^2 \,.
\end{equation}
Another preparatory step is a technical lemma:

\begin{lem}\label{Hilbert-coerc}
Let $(\mathbb{H},\sprod{\cdot}{\cdot})$ be a Hilbert space with induced norm $\|\cdot\|$. Let $v_1,v_2\in\mathbb{H}$ satisfy 
\begin{equation}\label{coerc-cond}
    1 + \sprod{v_1}{v_2}\ne 0 \,. 
\end{equation}    
Then there exists a constant $\kappa(\|v_1\|, \|v_2\|, \sprod{v_1}{v_2})>0$ such that
$$
   \|u + v_1\sprod{u}{v_2}\|^2 \ge \kappa \|u\|^2 \qquad\forall\, u\in \mathbb{H} \,.
$$ 
\end{lem}

\begin{proof}
It is easily seen that \eqref{coerc-cond} is needed for the definiteness of the left hand side, since it can only vanish if $u = \lambda v_1$, whence
$u + v_1\sprod{u}{v_2} = \lambda v_1 (1+ \sprod{v_1}{v_2})$. By homogeneity it suffices to show the inequality on the unit sphere. 
The left hand side can be written as
$$
   I(u) := \|u\|^2 + 2\sprod{u}{v_1}\sprod{u}{v_2} + |v_1|^2 \sprod{u}{v_2}^2 \,.
$$
With the decomposition $u = u_{12} + u^\bot$, $u_{12}\in{\rm sp}\{v_1,v_2\}$, $u^\bot\in{\rm sp}\{v_1,v_2\}^\bot$, we have
$$
  I(u) = I(u_{12}) + \|u^\bot\|^2 \,,
$$
showing that we may restrict our attention to $u\in{\rm sp}\{v_1,v_2\}\cap \{u\in\mathbb{H}:\, \|u\|=1\}$, a compact one- or two-dimensional manifold,
where the continuous functional $I$ assumes its minimum $\kappa\ge 0$, which cannot be equal to zero by definiteness. 
\end{proof}

\begin{lem}\label{lem:micro}{\bf (Microscopic coercivity)}
With the above definitions there exists $\lambda_m>0$ such that
\begin{equation}
- \sprod{\sfL F}{F} \ge \lambda_m \|(1-\Pi)F\|^2 \,.\label{est1}
\end{equation}
\end{lem}
\begin{proof}\label{proof mic coer}
By \eqref{lin-entr-diss} and \eqref{def LFF} we have
\begin{eqnarray}
   -\sprod{\sfL F}{F} &=& \|\hr M(\Tinf) - \hf + \hT g_1\|_{v,x}^2 + \frac{D}{\Tinf^2}\|\grad \hT\|_x^2 \,,\nn\\
   && \mbox{with } g_1(v) = \frac{\finf}{\Tinf}\left( \frac{|v|^2}{\Tinf} - \frac{d}{2}\right) \,.\label{lff}
\end{eqnarray} 
By the conservation laws \eqref{cons-zero} we have
\begin{equation}\label{cons-zero2}
   \frac{1}{|\Omega|} \ints{\hT} =: \bar T = -\sprod{f^\bot}{g_2} \quad\mbox{with } f^\bot := \hf - \hr M(\Tinf) \,,\, g_2 := \frac{|v|^2 \finf}{|\Omega|}  \,.
\end{equation}
The Poincar\'e inequality \eqref{Poincare} implies
\begin{eqnarray*}
   -\sprod{\sfL F}{F} &\ge& \left\|f^\bot + g_1\sprod{f^\bot}{g_2} - g_1(\hT - \bar T)\right\|_{v,x}^2 + \frac{Dc_P}{\Tinf^2}\left\|\hT-\bar T\right\|_x^2 \\
   &\ge& \left\|f^\bot + g_1\sprod{f^\bot}{g_2} \right\|_{v,x}^2 - 2\|g_1\|_v \left\|f^\bot + g_1\sprod{f^\bot}{g_2} \right\|_{v,x} \left\|\hT-\bar T\right\|_x  \\
   && + \left( \|g_1\|_v^2 + \frac{Dc_P}{\Tinf^2}\right) \left\|\hT-\bar T\right\|_x^2 \\
   &\ge& c_1 \left\|f^\bot + g_1\sprod{f^\bot}{g_2} \right\|_{v,x}^2 + c_1 \left\|\hT-\bar T\right\|_x^2 \,,
\end{eqnarray*}
for some $c_1>0$. A straightforward computation using \eqref{M-moments} gives $\sprod{g_1}{g_2}_{v,x} = \rinf d/2 \ne -1$, allowing the application of
Lemma \ref{Hilbert-coerc} to the first term on the \rhs with the consequence
$$
   -\sprod{\sfL F}{F} \ge c_1\kappa_1 \left\|f^\bot \right\|_{v,x}^2 + c_1 \left\|\hT-\bar T\right\|_x^2 \,.
$$
Since, by \eqref{cons-zero2}, we also have 
$$
   |\bar T| \le \left\| f^\bot \right\|_{v,x} \|g_2\|_v \,,
$$
we obtain
\begin{eqnarray*}
  -\sprod{\sfL F}{F} &\ge& \frac{c_1\kappa_1}{2} \left\|f^\bot \right\|_{v,x}^2 + c_1\left(\frac{\kappa_1}{2\|g_2\|_v^2} - |\Omega|\right)\bar T^2 
   + c_1 \|\hT\|_x^2 \\
   &=& \frac{c_1\kappa_1}{2} \left\|f^\bot \right\|_{v,x}^2 + c_1 \left\|\hT + \left( \sqrt{\frac{\kappa_1}{2|\Omega|\|g_2\|_v^2}} - 1\right)\bar T\right\|_x^2 \\
   &\ge& \frac{c_1\kappa_1}{2} \left\|f^\bot \right\|_{v,x}^2 + c_1 \kappa_2 \left\|\hT \right\|_x^2 \,,
\end{eqnarray*}
by another application of Lemma \ref{Hilbert-coerc}, now in $L^2(\Omega)$ with 
$$
   v_1= \sqrt{\frac{\kappa_1}{2|\Omega|\|g_2\|_v^2}} - 1 \,,\quad v_2 = \frac{1}{|\Omega|} \,,\quad 
     1 + \sprod{v_1}{v_2} = \sqrt{\frac{\kappa_1}{2|\Omega|\|g_2\|_v^2}} \ne 0 \,.
$$
We conclude with $\lambda_m = c_1 \min\{\kappa_1/2,\kappa_2\Tinf^2\}$, since \\$\|(1-\Pi)F\|^2 = \|f^\bot\|_{v,x}^2 + \|\hT\|_x^2/\Tinf^2$.
\end{proof}

\begin{remark}
In principle, the coercivity constant $\kappa$ in Lemma \ref{Hilbert-coerc} is computable and, thus, the same is true for $\lambda_m$.
However, we have not found a reasonably simple representation.
\end{remark}

The second term on the \rhs of \eqref{derivt entr} is expected to control $\Pi F$. Actually the operator 
$A\sfT\Pi = (1 + (\sfT\Pi)^*\sfT\Pi)^{-1}(\sfT\Pi)^*\sfT\Pi$ can be interpreted as application of the nonlinear function $z\mapsto \frac{z}{1+z}$ to the 
symmetric nonnegative operator $(\sfT\Pi)^*\sfT\Pi$. Thus, it suffices to show that the restriction of the latter to the null space of $\sfL$ has as a spectral
gap.

\begin{lem} \label{lem:macro}{\bf (Macroscopic coercivity)}
With the above definitions, we have
$$
  \|\sfT\Pi F\|^2 \ge \lambda_M \|\Pi F\|^2 \,,\qquad\mbox{with } \lambda_M = \frac{c_P\Tinf}{2}\,.
$$
\end{lem}

\begin{proof}
By $\Pi F = (\hr M(T_\infty), 0)$ we have $\sfT\Pi F = (v\cdot\nabla_x \hr M(T_\infty), 0)$ and, after a straightforward computation using \eqref{M-moments},
the Poincar\'e inequality \eqref{Poincare}, and \eqref{cons-zero},
$$
  \|\sfT\Pi F\|^2 = \frac{\Tinf}{2\rinf} \|\nabla_x \hr\|_x^2 \ge \frac{c_P\Tinf}{2\rinf} \|\hr\|_x^2 = \frac{c_P\Tinf}{2} \|\Pi F\|^2\,.
$$
\end{proof}

Considering the above discussion, it is an immediate corollary of Lemma \ref{lem:macro} that
\begin{equation}\label{macro-est}
   \langle \sfA\sfT\Pi F,F\rangle \ge \frac{\lambda_M}{1+\lambda_M} \|\Pi F\|^2 \,.
\end{equation}
Thus, the first two terms on the \rhs of \eqref{derivt entr} provide the desired coercivity. The remaining three terms are controlled by the first two for
$\delta$ small enough, if they can be estimated in terms of the product $\|F\|\,\|(1-\Pi)F\|$. For the third term this is not true in the general abstract
setting. It requires an additional algebraic property.

\begin{lem}\label{lem:diff-lim} {\bf (Diffusive macroscopic limit)}
With the above definitions, we have
$$
   \Pi\sfT\Pi = 0 \,.
$$
\end{lem}

\begin{proof}
This follows immediately from the expression for $\sfT\Pi$ in the previous proof, since the mean velocity of $M(\Tinf)$ vanishes.
\end{proof}

The macroscopic operator $\mathscr{L}:=(\sfT \Pi)^* \sfT \Pi = -\Pi\sfT^2\Pi$ acts only on $\cal{N}(\sfL)$, whence Lemma \ref{lem:diff-lim} implies
$$
   \sfT\sfA = -\sfT(1+\mathscr{L})^{-1} \Pi\sfT = (1-\Pi)\sfT\sfA(1-\Pi) \,.
$$
In \cite[Lemma 1]{DolMouSchmHypoMass} it has been shown that the operator norm of $\sfT\sfA$ is bounded by 1 and therefore
\begin{equation}\label{TA-bound}
  |\sprod{\sfT\sfA F}{F}| \le \|(1-\Pi)F\|^2 \,.
\end{equation}

It remains to find bounds for the last two terms in \eqref{derivt entr}. 

\begin{lem} \label{lem:aux-est} {\bf (Boundedness of auxiliary operators)}
With the above definitions, we have
$$
   |\sprod{\sfA\sfT(1-\Pi)F}{F}| +  |\sprod{\sfA\sfL F}{F}| \le C_M \|(1-\Pi)F\| \,\|\Pi F\|\,,
$$
with
$$
   C_M = \sqrt{d(d+2)\rinf}\, + \sqrt{d\rinf\max\left\{ 1, \frac{\rinf d}{2}\right\}} \,.
$$
\end{lem}

\begin{proof}
We start with the observation
\begin{eqnarray*}
  |\sprod{\sfA\sfT(1-\Pi)F}{F}| &= & |\sprod{(1-\Pi)F}{\sfT^2(1+\mathscr{L})^{-1} \Pi F}| \\
  &\le& \|(1-\Pi)F\| \left\|\sfT^2(1+\mathscr{L})^{-1} \Pi F\right\|  \,.
\end{eqnarray*}
Introducing $G:= (1+\mathscr{L})^{-1} \Pi F$, we have $G = \Pi G = (\rho_G M(\Tinf), 0)$ with 
\begin{equation}\label{G-equ}
   \rho_G - \frac{\Tinf}{2}\Delta_x \rho_G = \hr \,.
\end{equation}
Testing this against $\Delta_x \rho_G$ implies
\begin{equation}\label{G-est}
  \|\Delta_x \rho_G\|_x = \|\nabla_x^2 \rho_G\|_x \le \frac{2}{\Tinf} \|\hr\|_x = \frac{2\sqrt{\rinf}}{\Tinf} \|\Pi F\|\,.
\end{equation}
Finally,
\begin{eqnarray*}
  \|\sfT^2 G\|^2 &=& \intsv{((v\cdot\nabla_x)^2 \rho_G)^2 M(\Tinf)} 
  \le \intv{|v|^4 M(\Tinf)} \left\| \nabla_x^2 \rho_G\right\|_x^2 \\
  &=& \frac{d(d+2)\Tinf^2}{4} \left\| \nabla_x^2 \rho_G\right\|_x^2 \le d(d+2)\rinf \|\Pi F\|^2 \,,
\end{eqnarray*}
completing the proof of the first result.

For the second estimate we compute
\begin{eqnarray*}
  |\sprod{\sfA\sfL F}{F}| = \left| \sprod{(1+\mathscr{L})^{-1}\Pi\sfT\sfL F}{F}\right| = \left| \sprod{(\bfQ_L(F),0)}{TG}\right| \,.
\end{eqnarray*}
It is correct to replace the second component of $\sfL F$ by zero on the right hand side, since also the second component of $\sfT$ vanishes.
The first factor is estimated by 
\begin{eqnarray*}
  \|(\bfQ_L(F),0)\|^2 &=& \|\bfQ_L(F)\|_{v,x}^2 = \|\hf - \hr M(\Tinf)\|_{v,x}^2 + \frac{\rinf d}{2} \frac{\|\hT\|_x^2}{\Tinf^2} \\
  &\le& \max\left\{ 1, \frac{\rinf d}{2}\right\} \|(1-\Pi)F\|^2
\end{eqnarray*}
For the second factor we show similarly to above that 
$$
   \|\sfT G\|^2 \le \frac{d\Tinf}{2}  \|\nabla_x \rho_G\|_x^2 = - \frac{d\Tinf}{2} \sprod{\rho_G}{\Delta_x \rho_G}_x \le d \|\hr\|_x^2 
   = d\rinf \|\Pi F\|^2 \,,
$$
where we have used \eqref{G-est} and the obvious consequence $\|\rho_G\|_x \le \|\hr\|_x$ of \eqref{G-equ}, completing the proof.
\end{proof}

Using Lemma \ref{lem:micro}, \eqref{macro-est}, \eqref{TA-bound}, and Lemma \ref{lem:aux-est} in \eqref{derivt entr} we obtain
\begin{equation}
\label{hypocoercivity}
\derivt \scrH[F] \le -(\lambda_m - \delta) \| (1-\Pi)F \|^2 - \frac{\delta\lambda_M}{1+\lambda_M} \|\Pi F\|^2 + \delta C_M  \|(1-\Pi F)\| \|\Pi F\| \,.
\end{equation}
It is obvious that for $\delta$ small enough there exists $\kappa> 0$, such that
$$
  \frac{d}{dt} \scrH[F] \le - \kappa \|F\|^2 \le -\lambda \scrH[F]  \,,\qquad\mbox{with } \lambda = \frac{2\kappa}{1+\delta} \,,
$$
where in the last inequality \eqref{equival} has been used. An application of the Gronwall lemma and again of \eqref{equival} complete the proof of
Theorem \ref{Thm Hypo}.

\section{Proof of the global existence of the non-linear problem in one dimension}
\label{sec existence}

In this section we prove Theorem \ref{Thm_Ex}, the existence of global solutions in the one dimensional case. As a first step, the problem is regularized by cutting
off small temperatures in the Maxwellian:
$$
   M_\eps(T) := M(\theta_\eps(T)) \,,\quad \theta_\eps(T) := \max\{T,\eps\} \,,\qquad \eps>0 \,.
$$
For the regularized problem, a mild formulation of the Cauchy problem for the kinetic equation \eqref{Kin} reads
\begin{equation}\label{Kin-eps}
   f_\eps(x,v,t) = e^{-t}f_0(x-vt,v) + \int_0^t e^{s-t} \rho_\eps(x-v(t-s),s)M_\eps(T_\eps(x-v(t-s),s))(v)ds \,,
\end{equation}
and we recall the mild formulation of the heat conduction problem:
\begin{equation}\label{Heat-eps}
   T_\eps(x,t) = G(\cdot,t)*T_0 + \int_0^t G(\cdot,t-s)*\left(E_\eps - \frac{d}{2}\rho_\eps T_\eps\right)(\cdot,s)ds \,,
\end{equation}
where $\rho_\eps$ and $E_\eps$ are the zeroth and, respectively, second order velocity moments of $f_\eps$.

\begin{lem}\label{lem:exist-reg}
Let the assumptions of Theorem \ref{Thm_Ex} hold. Then for any $\eps>0$ the problem \eqref{Kin-eps}, \eqref{Heat-eps} has a unique solution 
$(f_\eps,T_\eps)$ with 
$$
   (1 + |v|^2 + \log f_\eps)f_\eps \in L^\infty((0,\tau); L^1(\mathbb{T}^1 \times \R)) \,,\qquad T_\eps \in C^{1/2}(\mathbb{T}^1 \times [0,\tau]) \,,
$$
uniformly in $\eps$.
\end{lem}

\begin{proof}
Using the fixed point form $f_\eps = F_\eps(f_\eps,T_\eps)$, $T_\eps = G(f_\eps,T_\eps)$ of \eqref{Kin-eps}, \eqref{Heat-eps}, we start with a local existence 
result by Picard iteration on the space $L^\infty([0,t_0]; L^1((1+|v|^2)dx\,dv)\times L^\infty((0,t_0)\times \mathbb{T}^1) =: L^\infty_t L^{1,2}_{x,v} \times L^\infty_{x,t}$ 
with $t_0$ small enough. This is a straightforward exercise. We only state the estimates
\begin{align*}
    \left\|F_\eps(f,T)\right\|_{L^\infty_t L^{1,2}_{x,v}}  \le \left\|f_0\right\|_{L^{1,2}_{x,v}} 
    + t_0 \left\|f\right\|_{L^\infty_t L^1_{x,v}}  \left(1+ \frac{d}{2} \|T\|_{L^\infty_{x,t}}\right)
\end{align*}
and, using Lemma \ref{lem:heat0},
$$
  \|G(f,T)\|_{L^\infty_{x,t}} \le \|T_0\|_{L^\infty_x} 
  + \sqrt{t_0}\, C_\tau \left\|f\right\|_{L^\infty_t L^{1,2}_{x,v}} \left(1+ \frac{d}{2} \|T\|_{L^\infty_{x,t}}\right) \,,
$$
showing that for $t_0$ small enough, a set of the form 
$$
   \left\{ (f,T):\, \left\|f\right\|_{L^\infty_t L^{1,2}_{x,v}} + \|T\|_{L^\infty_{x,t}} < R \right\}  \qquad\mbox{with}\quad
    R > \left\|f_0\right\|_{L^{1,2}_{x,v}} + \|T_0\|_{L^\infty_x} \,,    
$$
is mapped into itself by $(F_\eps,G)$. The contraction property of $(F_\eps,G)$ follows from similar estimates. Here the temperature cut-off is needed:
$$
    \left| M_\eps(T_1)(v) - M_\eps(T_2)(v)\right| \le \frac{c}{\eps^2}\left(1+|v|^2\right) M_\eps(\tilde T)(v) |T_1 - T_2|\,,
$$
with $\tilde T$ between $T_1$ and $T_2$.

The next step is to obtain estimates, which can be used for arbitrarily large times.
Obviously, the local solution satisfies $f_\eps\ge 0$ (thus $\rho_\eps,E_\eps \ge 0$) and therefore, by Corollary \ref{cor:heat}, $T_\eps\ge 0$. We also have
$$
   \frac{d}{dt} \int_{\mathbb{T}^1} \rho_\eps dx = 0 \,,\qquad \frac{d}{dt} \int_{\mathbb{T}^1} (E_\eps + T_\eps) dx 
   = \frac{d}{2}\int_{\mathbb{T}^1}\rho_\eps (\theta_\eps(T_\eps)-T_\eps)dx \le \frac{\eps d}{2} \|f_0\|_{L^1_{x,v}}\,,
$$
implying
\begin{equation} \label{feps-bound}
    \left\|f_\eps(\cdot,\cdot,t)\right\|_{L^{1,2}_{x,v}} \le \left( 1 + \frac{\eps\tau d}{2}\right)\left\|f_0\right\|_{L^{1,2}_{x,v}} \,.
\end{equation}
Therefore, since \eqref{Heat-eps} with $\rho_\eps=0$ provides an upper bound for $T_\eps$, by Lemma \ref{lem:heat0},
\begin{equation} \label{Teps-bound}
   \|T_\eps\|_{L^\infty_{x,t}} \le \|T_0\|_{L^\infty_x} + C_\tau \left\|f_0\right\|_{L^{1,2}_{x,v}} \,.
\end{equation}
The last two estimates prove that the solution can be extended to the whole time interval $(0,\tau)$, no matter how long it is. The H\"older continuity of
$T_\eps$ follows from Lemma \ref{lem:heat}. It remains to prove the $L\log L$ bound for $f_\eps$ by a entropy dissipation, which carries over to the 
regularized problem: 
\begin{align*}
\frac{d}{dt}\calH(f_\eps,T_\eps) = & - \intsv{ (f_\eps - \rho_\eps M(\theta_\eps(T_\eps))) \ln \frac{f_\eps}{\rho_\eps M(\theta_\eps(T_\eps))} } \\
& - D \ints{\frac{|\grad T|^2}{T^2}} + \ints{E_\eps\left(\frac{1}{\theta_\eps(T_\eps)} - \frac{1}{T_\eps}\right)} \le 0\,.
\end{align*}

\end{proof}

The proof of Theorem \ref{Thm_Ex} will be completed by passing to the limit $\eps\to 0$ in \eqref{Heat-eps} and in the distributional formulation of \eqref{Kin-eps}:
\begin{align}
  & -\int_{\mathbb{T}^1\times\R} f_0\phi(t=0)dx\,dv + \int_{\mathbb{T}^1\times\R\times (0,\tau)} f_\eps(\phi - \partial_t\phi - v\cdot\nabla_x\phi)dx\,dv\,dt \nonumber\\
  & = \int_{\mathbb{T}^1\times\R\times (0,\tau)} \rho_\eps(x,t) M(1)(w) \phi\left(x,w\sqrt{\theta_\eps(T_\eps)},t\right) dx\,dw\,dt \,,\label{Kin-distr}
\end{align}
with $\phi\in C^1_0(\mathbb{T}^1\times\R\times [0,\tau))$.

By Lemma \ref{lem:exist-reg} and the Arzel\`a-Ascoli theorem, uniform convergence of $T_\eps$ follows (when restricted to appropriate subsequences). 
On the other hand, multiplication of the kinetic equation by $|v|^q$ and integration gives
$$
   \frac{d}{dt} \int_{\mathbb{T}^1\times\R} |v|^q f_\eps dx\,dv = c_q \int_{\mathbb{T}^1} \rho_\eps\theta_\eps(T_\eps)^{q/2} dx 
     - \int_{\mathbb{T}^1\times\R} |v|^q f_\eps dx\,dv \,,
$$
with $c_q = \int_{\R} |w|^q M(1)(w)dw$. Since by \eqref{feps-bound} and \eqref{Teps-bound} the first term on the right hand side is uniformly bounded, the same is true for
$\int_{\mathbb{T}^1\times\R} |v|^q f_\eps dx\,dv$. This implies tightness of the sets $\{f_\eps:\, \eps>0\}$ and $\{|v|^2 f_\eps:\, \eps>0\}$. Furthermore the uniform
$L\log L$ bound gives uniform integrability of $\{f_\eps:\, \eps>0\}$ by the De La Vall\'ee-Poussin criterion and therefore weak relative compactness in $L^1_{x,v,t}$ 
by the Dunford-Pettis theorem (see, e.g., \cite[Theorems 3.1, 3.2]{FGLSR}). Let $f\in L^1_{x,v,t}$ denote an accumulation point of $\{f_\eps:\, \eps>0\}$. It satisfies the 
same bounds as $f_\eps$, i.e. $(1+|v|^q + \log f)f \in L^\infty_t L^1_{x,v}$. For the kinetic energy densities $E_\eps$ and $E$ of $f_\eps$ and, respectively, $f$ we have
(with a test function $\phi(x,t)$)
\begin{eqnarray*}
   \left| \int_0^\tau\int_{\mathbb{T}^1}\phi(E_\eps - E)dx\,dt\right| &\le& \int_0^\tau \int_{\mathbb{T}^1} |\phi| \int_{\R\setminus B_R} |v^2| (f_\eps + f)dv\,dx\,dt \\
   && + \left| \int_0^\tau \int_{\mathbb{T}^1} \phi \int_{B_R} |v^2| (f_\eps - f)dv\,dx\,dt \right|
\end{eqnarray*}
By the bound on the $q$th order moments, the first term on the right hand side is of the order $R^{2-q}$ uniformly in $\eps$. For fixed $R$, the second term tends 
to zero as $\eps\to 0$, if $f_\eps \rightharpoonup f$. This provides weak $L^1$ convergence of $E_\eps$ and, in the same way, for $\rho_\eps$.

The uniform convergence of $T_\eps$ and the weak $L^1$ convergence of $f_\eps$, $\rho_\eps$, and $E_\eps$ (for a suitable subsequence) are sufficient
for passing to the limit $\eps\to 0$ in \eqref{Heat-eps} and \eqref{Kin-distr}, completing the proof of Theorem \ref{Thm_Ex}.

\section*{Appendix -- the heat equation with $L^1$ data}

In this section, we prove nonnegativity and H\"older regularity for the solution of the initial value problem
\begin{equation}\label{heat-app}
    \partial_t T - D\partial_x^2 T = h -\rho T\,,\qquad T(t=0)=T_0 \,,
\end{equation}
on the one-dimensional torus, with $h\ge 0$, $h,\rho \in L_t^\infty L_x^1$ and continuously differentiable initial datum $T_0$.
We start with the case $\rho=0$.

\begin{lem}\label{lem:heat0}
Let $\tau>0$, $T_0\in L^\infty(\T^1)$, $h\in L^\infty((0,\tau);L^1(\T^1))$, and 
\begin{equation} \label{mild-heat-norho}
    T(x,t) :=  \int_0^1 G(x-y,t)T_0(y)dy + \int_0^t \int_0^1 G(x-y,t-s)h(y,s)dy\,ds \,,
\end{equation}
with Green's function $G$ given by \eqref{Green}. Then there exists a constant $C_\tau>0$ such that
$$
  \|T\|_{L^\infty((0,\tau)\times\T^1)} \le \|T_0\|_{L^\infty(\T^1)} + C_\tau \|h\|_{L^\infty((0,\tau);L^1(\T^1))} \,,
$$
with $C_\tau = O(\sqrt{\tau})$ as $\tau\to 0$.
\end{lem}

\begin{proof}
Since
$$
  \int_0^1 G(x-y,t)dy = 1 \qquad\mbox{and}\qquad |G(x,t)| \le \frac{c_\tau}{\sqrt{t}} \,, 
$$
we have
$$
   |T(x,t)| \le \|T_0\|_{L^\infty(\T^1)} + 2\sqrt{t}\,c_\tau \|h\|_{L^\infty((0,\tau);L^1(\T^1))} \,,
$$
completing the proof.
\end{proof}

\begin{lem}\label{lem:heat}
Let $\tau>0$, $T_0\in C^1(\T^1)$, $h\in L^\infty((0,\tau);L^1(\T^1))$, and let $T$ be given by \eqref{mild-heat-norho}.
Then there exists a constant $C_\tau>0$ such that
$$
  \|T\|_{C^{1/2}([0,\tau]\times\T^1)} \le C_\tau \left( \|T_0\|_{C^1(\T^1)} + \|h\|_{L^\infty((0,\tau);L^1(\T^1))} \right) \,.
$$
\end{lem}

\begin{proof}
With the periodic extensions of $T_0$ and $h$, we also have 
\begin{eqnarray*}
    T(x,t) &=&  \int_{\R} G_{\R}(x-y,t)T_0(y)dy + \int_0^t \int_{\R} G_{\R}(x-y,t-s)h(y,s)dy\,ds \\
    &=& (G_{\R}(\cdot,t)*T_0)(x) + \int_0^t (G_\R(\cdot,t-s)*h(\cdot,s))(x)ds =: T_1(x,t) + T_2(x,t)\,,
\end{eqnarray*}
with
$$
    G_{\R}(x,t) = (4D\pi t)^{-1/2} \exp\left(- \frac{x^2}{4Dt}\right) \,.
$$
By 
$$
    \partial_x T_1 = G_{\R}*T_0'
$$
and the maximum principle, the derivative of $T_1$ with respect to $x$ is bounded. 
On the other hand, using
$$
   \erfc(x) = \frac{2}{\sqrt{\pi}} \int_{x}^{+\infty} e^{-s^2}\ds \qquad\mbox{with}\quad
\frac{d}{dx} \erfc\bigg(\frac{x}{\sqrt{4Dt}}\bigg) = - 2 G_\mathbb{R}(x,t)\,,
$$
integrations by parts with respect to $y$ on $[0,\infty)$ and on $(-\infty,0]$ give
$$
    T_1 (x,t) = T_0(x) - \frac{1}{2} \int_0^\infty \erfc\left(\frac{y}{\sqrt{4Dt}}\right) \left(T_0'(x-y) - T_0'(x+y)\right)dy
$$
This implies, with $t_2>t_1$,
\begin{eqnarray*}
   |T_1(x,t_1)-T_1(x,t_2)| &\le& c_1 \int_0^\infty \left( \erfc\left(\frac{y}{\sqrt{4Dt_2}}\right) - \erfc\left(\frac{y}{\sqrt{4Dt_1}}\right) \right)dy \\
   &=& c_2\left( \sqrt{t_2} - \sqrt{t_1}\right) \le c_2 \sqrt{t_2-t_1}\,,
\end{eqnarray*}
proving the claim for $T_1$.

For the $x$-dependence of $T_2$ we have
\begin{align*}
  & |T_2(x_1,t) - T_2(x_2,t)| \le \int_0^t \int_0^1 \left| G(x_1-y,s) - G(x_2-y,s)\right| |h(y,t-s)|dy\,ds \\
   & \hskip 2cm \le  \int_0^t \sup_{y\in (0,1)} \left| G(x_1-y,s) - G(x_2-y,s)\right| ds \, \|h\|_{L^\infty((0,\tau);L^1(\T^1))}
\end{align*}
H\"older regularity will follow from (with $0<\alpha<1$, $z_1+k=\xi_k\sqrt{s}$, $z_2+k=\eta_k\sqrt{s}$)
\begin{eqnarray*}
 && \frac{\left| G(z_1,s) - G(z_2,s)\right|}{|z_1-z_2|^\alpha} \le (4D\pi s)^{-1/2} 
      \sum_{k\in \mathbb{Z}} \frac{| e^{-(z_1+k)^2/(4D s)} - e^{-(z_2+k)^2/(4D s)}|}{|z_1-z_2|^\alpha} \\
  && = \frac{1}{\sqrt{4D\pi}} s^{-1/2-\alpha/2} 
     \sum_{k\in \mathbb{Z}} \frac{| e^{-\xi_k^2/8} - e^{-\eta_k^2/8}|}{|\xi_k-\eta_k|^\alpha}\left(e^{-(z_1+k)^2/(8s)} + e^{-(z_2+k)^2/(8s)}\right) \\ 
   && \le \frac{c_\alpha C_\tau}{\sqrt{4D\pi}} s^{-1/2-\alpha/2} \,,
\end{eqnarray*}
where $c_\alpha$ is the H\"older constant of $\xi\mapsto e^{-\xi^2/8}$ and $C_\tau$ is an estimate for the sum which remains. Since the right hand side is integrable
with respect to $s$ for $\alpha<1$, it follows that $T_2$ is H\"older continuous as a function of $x$ with any exponent $\alpha<1$.

Similarly, for the $t$-dependence of $T_2$, with $t_2>t_1$,
\begin{align}
  & |T_2(x,t_1) - T_2(x,t_2)| \le \int_0^{t_1} \int_0^1 \left| G(y,t_1-s) - G(y,t_2-s)\right| |h(x-y,s)|dy\,ds \nonumber\\
  & \hskip 3.5cm +  \int_{t_1}^{t_2} \int_0^1 G(y,t_2-s) |h(x-y,s)|dy\,ds \nonumber\\
   & \hskip 2cm \le  \int_0^{t_1} \sup_{y\in (0,1)} \left| G(y,t_1-s) - G(y,t_2-s)\right| ds \, \|h\|_{L^\infty((0,\tau);L^1(\T^1))} \nonumber\\
   & \hskip 3.5cm +  \int_{t_1}^{t_2} \sup_{y\in (0,1)} G(y,t_2-s)ds  \, \|h\|_{L^\infty((0,\tau);L^1(\T^1))} \,. \label{long-est}
\end{align}
We start with the easier last term: It is easily seen that $G(x,t)\le c(1+t^{-1/2})$ and, thus,
$$
   \int_{t_1}^{t_2} \sup_{y\in (0,1)} G(y,t_2-s)ds \le C_\tau \sqrt{t_2-t_1} \,.
$$
With the abbreviations $z_k := (x+k)^2/(4D)$, $t_j-s =: u_j^{-2}$, we have
$$
   \left| G(y,t_1-s) - G(y,t_2-s)\right| \le c \sum_{k\in\mathbb{Z}} \left| u_1e^{-u_1^2 z_k} - u_2 e^{-u_2^2 z_k}\right| \,.
$$
For the function $h_k(u) := ue^{-u^2 z_k}$ we have
$$
   |h_k'(u)| = \left|\left(1-2u^2 z_k\right)e^{-u^2 z_k}\right| \le c\, e^{-u^2 z_k/2} \,,
$$
and therefore
\begin{eqnarray*}
  && \left| G(y,t_1-s) - G(y,t_2-s)\right| \le c(u_1-u_2) \sum_{k\in\mathbb{Z}} \exp\left(-\tilde u z_k/2\right) \\
  && = c \left( \frac{1}{\sqrt{t_1-s}} - \frac{1}{\sqrt{t_2-s}}\right) \sqrt{\tilde t}\, G(y,2\tilde t) \le C_\tau \left( \frac{1}{\sqrt{t_1-s}} - \frac{1}{\sqrt{t_2-s}}\right) \,,
\end{eqnarray*}
with $\tilde t$ between $t_1-s$ and $t_2-s$. Integration with respect to $s$ shows that the first term on the right hand side of \eqref{long-est} can be estimated
by $c\sqrt{t_2-t_1}$, completing the proof.
\end{proof}

\begin{cor}\label{cor:heat}
Let $\tau>0$, $T_0\in C^1_+(\T^1)$, $h,\rho\in L^\infty((0,\tau);L^1(\T^1))$, and $h\ge 0$.
Then \eqref{heat-app} has a unique, nonnegative, mild solution $T\in C^{1/2}([0,\tau]\times\T^1)$,
satisfying
\begin{equation}\label{mild-heat-app}
    T(x,t) :=  \int_0^1 G(x-y,t)T_0(y)dy + \int_0^t \int_0^1 G(x-y,t-s)(h-\rho T)(y,s))dy\,ds \,,
\end{equation}
with Green's function $G$ given by \eqref{Green}. 
\end{cor}

\begin{proof}
For regularized versions of $h$ and $\rho$, standard results for parabolic equations can be applied,
giving a smooth, nonnegative solution satisfying \eqref{mild-heat}. With $\rho=0$ we obtain an upper solution,
which is uniformly bounded by Lemma \ref{lem:heat}. Another application of Lemma \ref{lem:heat} to \eqref{mild-heat-app}
with regularized data provides a uniform $C^{1/2}$ estimate for $T$, implying
uniform convergence in the limit, when the regularization is removed. This is sufficient for passing to the limit in  $h-\rho T\in L^\infty((0,\tau);L^1(\T^1))$.
\end{proof}
\bigskip

\noindent {\bf Acknowledgments.}~~This work has been supported by the Austrian Science Fund, grants no. W1245 and F65, and by the Humboldt foundation. 
G.F.~also thanks the Vienna School of Mathematics.


\pagestyle{plain}
\pagenumbering{roman}

\bibliographystyle{acm}

\end{document}